\documentclass[11pt]{article}

\usepackage{amsfonts,amssymb}
\usepackage{graphicx}

\newcommand{\N}{\mathbb{N}}
\newcommand{\C}{\mathbb{C}}
\newcommand{\Z}{\mathbb{Z}}
\newcommand{\R}{\mathbb{R}}

\newcommand{\Hom}{\mbox{Hom}}
\newcommand{\id}{\mbox{id}}

\newtheorem{theorem}{Theorem}[section]

\newtheorem{lemma}[theorem]{Lemma}
\newtheorem{proposition}[theorem]{Proposition}

\newtheorem{definition}[theorem]{Definition}

\newtheorem{example}[theorem]{Example}

\newtheorem{remark}[theorem]{Remark}

\input{epsf.sty}

\def\ch{ {\rm char }}

\def\id{ {|} }

\begin{document}

\title{Cohomology of Frobenius Algebras and the Yang-Baxter Equation}

\author{
J. Scott Carter\footnote{Supported in part by
NSF Grant DMS  \#0603926.}
\\ University of South Alabama
\and
Alissa S. Crans \\
Loyola Marymount University
\and
Mohamed
Elhamdadi
\\ University of South Florida
\and
Enver Karadayi
\\ University of South Florida
\and
Masahico Saito\footnote{Supported in part by NSF Grant DMS
\#0603876.}
\\ University of South Florida
}

\maketitle
\begin{center}
{\it Dedicated to the memory of Xiao-Song Lin}
\end{center}

\begin{abstract}
A cohomology theory for multiplications and comultiplications of Frobenius
algebras is developed
in 
low dimensions
in analogy with Hochschild cohomology of
bialgebras 
based on
deformation theory.
Concrete computations are provided for key examples.

Skein theoretic constructions give rise to solutions to the
Yang-Baxter equation
using multiplications and comultiplications of Frobenius
algebras, and $2$-cocycles are used to obtain deformations
of $R$-matrices thus obtained.
\end{abstract}

\section{Introduction}

Frobenius algebras are interesting to topologists as well as
algebraists
for 
numerous reasons including the following.
First,
$2$-dimensional topological quantum field theories are formulated
in terms of commutative Frobenius 
algebras
(see \cite{Kock}).
Second, a Frobenius algebra structure exists on any
finite-dimensional Hopf algebra with a left integral defined in the dual space.
These
Hopf algebras have found applications in topology through
Kuperberg's invariant
\cite{Kup91,Kup96},
the Henning invariant   \cite{KR95,Oht},
and the theory of quantum groups from which the post-Jones invariants arise.
Third, there is a $2$-dimensional Frobenius algebra that underlies Khovanov's  cohomology theory~\cite{Kh}. See also~\cite{AF}.

Our interest herein is to extend the cohomology theories defined in \cite{CCES1,CCES2} to Frobenius algebras and thereby construct new solutions to the Yang-Baxter equation (YBE).  We expect that there are connections among these cohomology theories that extend beyond their formal definitions. Furthermore, we anticipate topological, categorical, and/or physical applications because of the diagrammatic nature of the theory.

The $2$-cocycle conditions of Hochschild cohomology of algebras
and
bialgebras can be  interpreted
via deformations of algebras~\cite{GS}.
In other words, a map
satisfying the associativity condition can be deformed to obtain a
new associative map in a larger vector space using $2$-cocycles.
The same interpretation can be applied to quandle cohomology
theory~\cite{CJKLS,CENS,CES}.
A quandle is a set
equipped with  a self-distributive binary
operation
satisfying a few additional conditions
that correspond to the properties 
that conjugation in a group enjoys. 
Quandles have been used in
knot theory extensively (see \cite{CJKLS} and references therein
for  more aspects of quandles). Quandles and related structures
can be used to
construct 
set-theoretic solutions
(called  $R$-matrices)
to the Yang-Baxter
equation (see, for example, \cite{K&P} and its references). {}From
this point of view, combined with the deformation $2$-cocycle
interpretation, a quandle $2$-cocycle can be regarded as giving a
cocycle deformation of an $R$-matrix. Thus we
extend this idea to other algebraic constructions of
$R$-matrices
and  construct new $R$-matrices from old via
$2$-cocycle deformations.

In \cite{CCES1,CCES2}, new $R$-matrices were constructed via $2$-cocycle deformations in two other algebraic contexts.
Specifically,
in \cite{CCES1}, self-distributivity was revisited from
the point of
view of coalgebra categories, thereby unifying Lie algebras and
quandles in these categories. Cohomology theories of Lie algebras
and quandles were given
via a single
definition,
and deformations of
$R$-matrices were constructed. In \cite{CCES2}, the adjoint map of
Hopf algebras,
which
corresponds to the group conjugation  map,
was studied from the same viewpoint.  
A cohomology theory was
constructed based on equalities satisfied by the adjoint map that are sufficient for
it to satisfy the YBE.

In this paper, we present an analog for Frobenius algebras according to the following organization. After a brief review of necessary materials in Section~\ref{prelimsec},
a cohomology theory for Frobenius algebras is constructed in Section~\ref{cohsec}
via deformation theory.
Then Yang-Baxter solutions
are constructed by skein methods in Section~\ref{YBEsec},
followed by 
deformations of $R$-matrices by
$2$-cocycles.

The reader should be aware
that the composition of the maps is
read in the standard way from right to left $(gf)(x)=g(f(x))$ in
text and from bottom to top in the diagrams. In this way, when
reading from left to right one can draw from top to bottom and
when reading a diagram from top to bottom, one can display the
maps from left to right. The argument of a function (or input
object from a category)  is found at the bottom of  the diagram.

\section{Preliminaries} \label{prelimsec}

A {\it Frobenius algebra} is an (associative) algebra (with
multiplication $\mu: A \otimes A \rightarrow A$ and unit $\eta: k
\rightarrow A$) over a field $k$ with a nondegenerate associative
pairing $\beta:  A \otimes A \rightarrow k$. Throughout
this paper all algebras are
finite-dimensional unless
specifically
stated otherwise. The  pairing $\beta$ is also expressed by
$\langle x|y\rangle=\beta(x \otimes y)$ for $x, y \in A$, and it
is {\it associative}
in the sense that
$\langle xy | z\rangle=\langle x|y
z\rangle$ for any $x,y,z \in A$.

\begin{figure}[htb]
\begin{center}
\includegraphics[width=4.5in]{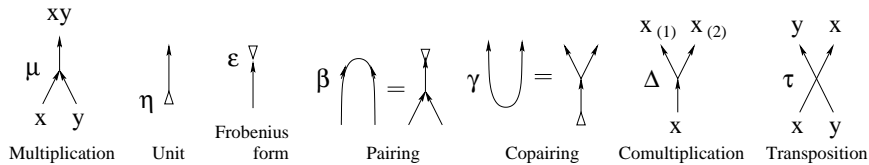}
\end{center}
\caption{Diagrams for Frobenius algebra maps}
\label{frodiag}
\end{figure}

A Frobenius algebra $A$ has a linear functional $\epsilon: A
\rightarrow k$,
called the {\it Frobenius form},
such that the kernel
contains no nontrivial  left ideal.
It is defined from $\beta$ by $\epsilon(x)=\beta(x \otimes 1)$,
and conversely, a Frobenius form gives rise to a nondegenerate
associative pairing $\beta$ by $\beta(x \otimes y)=\epsilon(xy)$,
for $x, y \in A$.

\begin{figure}[htb]
\begin{center}
\includegraphics[width=4in]{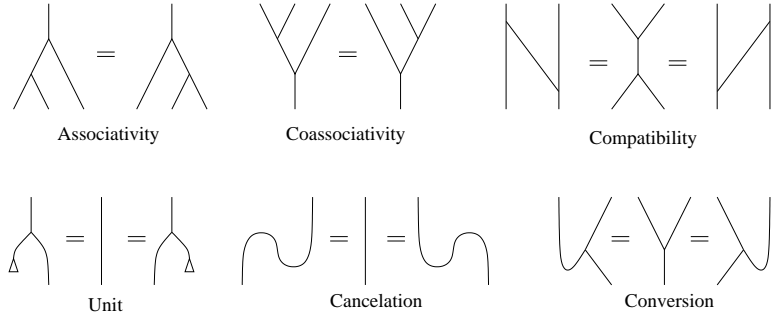}
\end{center}
\caption{Equalities among Frobenius algebra maps }
\label{froaxioms}
\end{figure}

A Frobenius form has a unique copairing
$\gamma: k \rightarrow A \otimes A$
characterized by
$$ (\beta \otimes |)(| \otimes \gamma) = |
= (| \otimes \beta)(\gamma \otimes |) , $$
where $|$ denotes the identity homomorphism on the algebra.  We call this relation the
{\it cancelation} of $\beta $ and $\gamma$.
 See the middle entry in the bottom row of Fig.~\ref{froaxioms}.
This notation will distinguish this function from the identity element $1=1_A=\eta(1_k)$ 
of the algebra that is the image of the identity of the ground field.
A Frobenius algebra $A$ determines a coalgebra structure with
$A$-linear (coassociative) comultiplication and the counit
defined using the Frobenius form. The comultiplication  $\Delta: A
\rightarrow A \otimes A$ is defined by
\begin{eqnarray*}
\Delta  &=&  (\mu \otimes |)(| \otimes \gamma) \\
&=& (| \otimes \mu)(\gamma \otimes |).
\end{eqnarray*}
The multiplication and comultiplication
satisfy the following equality:
$$ \Delta \mu = (\mu \otimes |)(| \otimes \Delta )=(| \otimes \mu)(\Delta \otimes |)$$
which we call the {\it Frobenius compatibility condition}.

A Frobenius algebra is {\it symmetric} if the pairing is
symmetric,
meaning that $\beta(x \otimes y)=\beta(y \otimes x)$ for any $x, y
\in A$. A Frobenius algebra is {\it commutative} if it is commutative as
an algebra.   It is known (\cite{Kock} Prop. 2.3.29) that
a Frobenius algebra is commutative if and only if it is
cocommutative as a coalgebra.

The map $\mu \Delta $ of a Frobenius algebra
is called the
{\it handle operator,} and
corresponds to multiplication by a central element called the {\it handle element}
$\delta_h =\mu \gamma (1) $ (\cite{Kock}, page 128).

Any semisimple Hopf algebra gives rise to
a Frobenius algebra structure (see, for example, \cite{Kock}, page $135$).
Let $H$ be a finite-dimensional Hopf algebra
with multiplication $\mu$ and unit $\eta$.
Then $H$ is semisimple if and only if the bilinear form
$\beta_{a\; b}=\sum  \mu_{c\; d}^d \ \mu_{a\; b}^c $ is nondegenerate.
If $H$ is semisimple, then the above defined $\beta$ gives rise to a Frobenius pairing.
In this case the Frobenius form (a counit of the Frobenius algebra structure)
is defined by $\epsilon_a=\sum  \mu_{c\; d}^d =T$,
the trace of $H$.
This counit and the induced comultiplication of the resulting Frobenius algebra
structure should not be confused with the counit and the comultiplication
of the original Hopf algebra.  
We thank Y. Sommerhausser for explaining these relationships to us.

We compute cohomology groups and the Yang-Baxter solutions
for a variety of examples, 
and we review these
mostly from \cite{Kock}.
{}From the point of view of TQFTs, the value $\epsilon \eta (1)$ corresponds to a sphere
$S^2$,  $\delta_0 = \beta \gamma (1)$ to a torus $T^2$,
and $\mu \Delta $ to adding an extra $1$-handle to a tube, so we
will compute these values and maps.

\begin{example}{\bf Complex numbers with trigonometric comultiplication.}\label{cplxex}
{\rm
Let $A=\C$ over $k=\R$ and
let the basis be denoted by $1$ and $i=\sqrt{-1}$. Then the
Frobenius form $\epsilon$ defined by
$\epsilon(1)=1$ and
$\epsilon(i)=0$ gives rise to the comultiplication $\Delta$,
which is Sweedler's trigonometric coalgebra
with
 $\Delta(1)=1 \otimes 1 - i \otimes i $, and
$\Delta(i)=i \otimes 1 + 1 \otimes i $.
We compute  $\epsilon \eta (1) = 1$,
$\mu \Delta  (1)=2 \eta(1) $,
and $\mu \Delta (i)=2i \eta(1) $
so that $\mu \Delta $ is
multiplication by $2$, which is the handle element $\delta_h$.
We also have $\delta_0=2$. } \end{example}

\begin{example}{\bf Polynomial algebras.\/}\label{polyex} {\rm
Polynomial rings $k[x]/(x^n)$ over a field $k$ where $n$ is a positive
integer are Frobenius
algebras. In particular, for $n=2$, the algebra $A=k[x]/(x^2)$ was
used in the Khovanov homology of knots~\cite{Kh}. For
$A=k[x]/(x^2)$,
the Frobenius form $\epsilon: A \rightarrow k$ is
defined by $\epsilon(x)=1$ and $\epsilon(1)=0$.
This induces the comultiplication
$\Delta: A \rightarrow A \otimes A$ determined by
$\Delta(1)=1 \otimes x + x \otimes 1$ and
$\Delta(x)= x \otimes x $. The handle element is $\delta_h=2x$.

More generally for $A=k[x]/(x^n)$
and $\epsilon(x^{j})=1$ for $j=n-1$ and $0$ otherwise,
the
comultiplication is determined by $\Delta(1) =
\sum_{i=0}^{n-1} x^i \otimes x^{n-1-i}$.
We have
$\mu \Delta (1) = n x^{n-1} $ and
the handle element is $\delta_h=n x^{n-1} $.

} \end{example}

\begin{example}{\bf Group algebras.\/}\label{gpex} {\rm
The group algebra $A=kG$ for a finite group $G$ over a field $k$
is
a Frobenius algebra with $\epsilon (x)=0$ for any $G \ni x \neq 1$
and $\epsilon (1)=1$, where $1$ is identified with the identity
element. The induced comultiplication is given by $\Delta(x)=\sum_{yz=x} y
\otimes z$.

One computes $\epsilon \eta (1) = 1$ and $\mu \Delta (x)=|G|x $,
where $|G|$ is the order of $G$. In particular, note that $\mu
\Delta = \delta_1 \id$ (recall that $\id$ denotes the identity map), where $\delta_1=|G|$ (the order of the group $G$), and $(\mu \Delta )^n =
\delta_1^n  \id$ for any $n \in \N$,
so that the handle element is $\delta_h=\delta_1=|G|$, and
$\delta_0=|G|$.

There are other Frobenius forms
on the group algebra A (again from \cite{Kock}). For example, for
$A=kG$ where $G$ is the symmetric group on three letters,
$A=k\langle x, y\rangle / (x^2-1, y^2-1, xyx=yxy )$, and
$\epsilon(xyx)=1$ and otherwise zero, is a Frobenius form and the
handle element is $2(xyx+x+y)$. } \end{example}

\begin{example}{\bf $q$-Commutative polynomials.\/}\label{qpolyex} {\rm
Let $X= k \langle x, y \rangle /(x^2, y^2, yx -  q xy) $
where $q=- A^{-2}$ for $A \in k$ with
polynomial multiplication and $\epsilon(xy)=iA$,
and  zero for other basis elements.
Then
$$ \gamma(1)( =\Delta\eta (1)) = iA (x \otimes y) - i A^{-1} (y \otimes x)
- i A^{-1} (xy \otimes 1 + 1 \otimes xy ) .$$
One computes
\begin{eqnarray*}
\Delta(x) &=& - iA^{-1} (x \otimes xy + xy \otimes x) , \\
\Delta(y) &=& - iA^{-1} (y \otimes xy + xy \otimes y) , \\
\Delta(xy) &=& - iA^{-1} (xy \otimes xy) .
\end{eqnarray*}
The handle element is $\delta_h = iA^{-1} (A - A^{-1} )^2 (xy) $.
} \end{example}

\section{Deformations and cohomology groups}\label{cohsec}

We describe the deformation theory of multiplication and
comultiplication for Frobenius algebras mimicking \cite{GS},
\cite{MrSt}, and our approach in \cite{CCES1,CCES2}. This approach
will yield the definition of $2$-cocycles. We will define the
chain complex for Frobenius algebras with chain groups in low
dimensions \cite{GS, MrSt}. We expect topological applications in low dimensions.
The differentials are defined via diagrammatically defined
identities among relations.

\subsection{Deformations}\label{deformsec}

In \cite{MrSt},
deformations of  bialgebras were described.
We follow that formalism and
give deformations of multiplications and comultiplications of Frobenius algebras.
A deformation of $A=(V, \mu, \Delta) $ is a
$k[[t]]$-Frobenius algebra $A_t=(V_t, \mu_t, \Delta_t)$, where
$V_t=V \otimes k[[ t ]]$ and
$ V_t/(tV_t) \cong V$.
Deformations
of $\mu$ and  $\Delta$ are given by $\mu_t= \mu + t \mu_1 + \cdots
+ t^n \mu_n + \cdots : V_t \otimes V_t \rightarrow V_t$ and
$\Delta_t = \Delta + t \Delta_1 + \cdots + t^n \Delta_n + \cdots :
V_t \rightarrow V_t \otimes V_t$ where $\mu_i: V \otimes V
\rightarrow V $, $\Delta_i : V \rightarrow V \otimes V$, $i=1, 2,
\cdots$, are sequences of maps. Suppose $\bar{\mu}=\mu + \cdots +
t^n \mu_n$ and $\bar{\Delta} =\Delta + \cdots + t^n \Delta_n$
satisfy the
Frobenius
conditions (associativity, compatibility,
and coassociativity) mod $t^{n+1}$, and suppose that  there exist
$\mu_{n+1}: V \otimes V \rightarrow V$ and $\Delta_{n+1}: V
\rightarrow V \otimes V$ such that $\bar{\mu}+t^{n+1} \mu_{n+1}$
and $\bar{\Delta}+ t^{n+1} \Delta_{n+1}$ satisfy the Frobenius
algebra conditions mod $t^{n+2}$. Define  $\xi_1 \in
\Hom(V^{\otimes 3}, V)$, $\xi_2, \xi_2' \in \Hom(V^{\otimes 2}, V^{\otimes
2})$, and $\xi_3 \in \Hom(V, V^{\otimes 3})$ by:
\begin{eqnarray*}
\bar{\mu} (\bar{\mu} \otimes |) - \bar{\mu}(| \otimes \bar{\mu})
&=& t^{n+1}\xi_1
\quad {\rm mod}\  t^{n+2} , \label{hoch2d1} \\
\bar{\Delta} \bar{\mu} -
(\bar{\mu} \otimes |)(| \otimes \bar{\Delta })&=&  t^{n+1} \xi_2
\quad {\rm mod}\  t^{n+2} ,  \\
\bar{\Delta} \bar{\mu} -
( | \otimes \bar{\mu})( \bar{\Delta }  \otimes | )
&=&  t^{n+1} \xi_2'
\quad {\rm mod}\  t^{n+2} , \label{hoch2d2}\\
(\bar{\Delta} \otimes |) \bar{\Delta} - (| \otimes \bar{\Delta})\bar{\Delta}
&=&  t^{n+1} \xi_3
\quad {\rm mod}\  t^{n+2} . \label{hoch2d3}
\end{eqnarray*}

\begin{remark}\label{deformrem}
{\rm
The operators in the quadruple $(\xi_1, \xi_2, \xi_2', \xi_3)$ 
form the {\it primary obstructions}
to formal deformations of multiplication and comultiplication
of a Frobenius algebra \cite{MrSt}.
} \end{remark}

For the associativity of $\bar{\mu}+t^{n+1} \mu_{n+1}$ mod
$t^{n+2}$ we obtain:
$$(\bar{\mu}+t^{n+1} \mu_{n+1})((\bar{\mu}+t^{n+1} \mu_{n+1}) \otimes |)-
(\bar{\mu}+t^{n+1} \mu_{n+1})(| \otimes (\bar{\mu}+t^{n+1}
\mu_{n+1}))=0 \  {\rm mod}\  t^{n+2} $$ which is equivalent by
degree calculations to:
\begin{eqnarray}
(d^{2,1}(\mu_{n+1})= ) & &   \mu (| \otimes \mu_{n+1})
+ \mu_{n+1}(| \otimes \mu )
  - \mu (\mu_{n+1} \otimes |)
  -\mu_{n+1}(\mu \otimes |)
=\xi_1,
\label{fro2cocy1}
\end{eqnarray}
where $d^{2,1}$ is one of the differentials we will define in the following section.
Similarly, from the Frobenius compatibility condition and
coassociativity we obtain
\begin{eqnarray}
(d^{2,2}_{(1)}(\mu_{n+1}, \Delta_{n+1})= ) & &
\Delta \mu_{n+1}
+ \Delta_{n+1} \mu -
(\mu \otimes |)(| \otimes \Delta_{n+1} )
- (\mu_{n+1} \otimes |)(| \otimes \Delta )=\xi_2,
\label{fro2cocy21}\\
(d^{2,2}_{(2)}(\mu_{n+1}, \Delta_{n+1})= ) & &
\Delta \mu_{n+1}
+ \Delta_{n+1} \mu -
(| \otimes \mu)(\Delta_{n+1}  \otimes |)
- ( | \otimes \mu_{n+1})(  \Delta \otimes | )=\xi_2',
\label{fro2cocy22}\\
(d^{2,3}(\Delta_{n+1})= ) & &
(\Delta \otimes |) \Delta_{n+1}
+  (\Delta_{n+1} \otimes |) \Delta
- (| \otimes \Delta) \Delta_{n+1}
- (| \otimes \Delta_{n+1} ) \Delta=
\xi_3,
\label{fro2cocy3}
 \end{eqnarray}
 where there are two types of compatibility conditions
 for $d^{2,2}$.

In summary we proved the following:
\begin{lemma} \label{deformlem}
The maps   $\bar{\mu}+t^{n+1} \mu_{n+1}$
and $\bar{\Delta}+t^{n+1} \Delta_{n+1}$  satisfy the
associativity, coassociativity and
Frobenius compatibility  conditions mod $t^{n+2}$ if and only if
the equalities (\ref{fro2cocy1}), (\ref{fro2cocy21}), (\ref{fro2cocy22}) and
(\ref{fro2cocy3}) are satisfied.
\end{lemma}

\subsection{Chain groups}

Let $A$ be a Frobenius algebra.
We define chain groups as follows.  
\begin{eqnarray*}
C^{n, i}_f (A;A) & =&  \Hom(A^{\otimes (n+1-i)}, A^{\otimes i} ), \\
C^n_f (A;A) & =& \oplus_{ 0 <
i \leq n} \ C^{n, i}_f (A;A) .
\end{eqnarray*}
Specifically, chain groups in low dimensions of our concern are:
\begin{eqnarray*}
C^1_f (A;A) &=& \Hom(A, A), \\
C^2_f (A;A) &=& \Hom(A^{\otimes 2}, A)\oplus \Hom(A, A^{\otimes 2}), \\
C^3_f (A;A) &=& \Hom(A^{\otimes 3}, A)\oplus
\Hom(A^{\otimes 2}, A^{\otimes 2})\oplus \Hom(A, A^{\otimes 3}).
\end{eqnarray*}
In the remaining sections we will define differentials that are
homomorphisms between the chain groups:
$$d^{n, i}_f =  d^{n, i}:  C^n_f (A;A)
\rightarrow C^{n+1, i}_f (A;A) (= \Hom(A^{\otimes (n+2-i)}, A^{\otimes i} ) )$$
which  will be  defined individually for $n=1,2,3$ and for $i$ with $0 \leq i \leq n$,
and
\begin{eqnarray*}
D_1 & =& d^{1, 1}-d^{1,2} : C^1_f (A;A)  \rightarrow C^2_f (A;A) ,\\
D_2^{(i)}&=& d^{2, 1} + d^{2, 2}_{(i)} + d^{2, 3}:  C^2_f (A;A)  \rightarrow C^3_f (A;A)
,\\
D_3&=& d^{3, 1} +d^{3, 2} +d^{3, 3} +d^{3, 4}:
 C^3_f (A;A)  \rightarrow C^4_f (A;A) .
\end{eqnarray*}
Define $C^0_f(A;A)=0$ by convention. {}From now on
the subscripts
$f$ for differentials are omitted for simplicity if no confusion
arises.

\subsection{First differentials}

By analogy with the differential for associative multiplication, we
make the following definition:

\begin{figure}[htb]
\begin{center}
\includegraphics[width=5in]{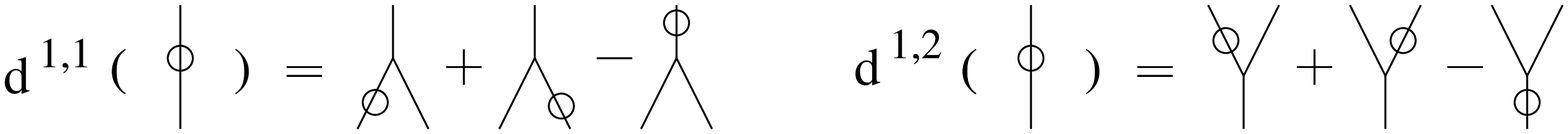}
\end{center}
\caption{First differentials }
\label{d1}
\end{figure}

\begin{definition}{\rm
The first differentials $d^{1,1}:  C^{1, 1}_f(A;A)   \rightarrow
C^{2,1}_f (A;A)  $ and $d^{1,2}:  C^{1, 1}_f(A;A)  \rightarrow
C^{2,2}_f (A;A)  $
are defined, respectively,  by
\begin{eqnarray*}
d^{1,1}(h) &=&  \mu (h \otimes |)+ \mu (| \otimes h) - h\mu  , \\
d^{1,2}(h) &=&  (h \otimes |)\Delta +  (| \otimes h) \Delta-
\Delta h .
\end{eqnarray*}
Then define $D_1: C^1_f (A;A)  \rightarrow C^2_{f} (A;A)$ by
$D_1=d^{1,1}-d^{1,2}$.
} \end{definition}
Diagrammatically, we represent $d^{1,i}$ for $i=1,2$, as depicted in
Fig.~\ref{d1}. A map $h \in C^{1, 1}_f(A;A)$ is represented by a
white circle on a vertical string and the multiplication and comultiplication
are depicted by
two distinct trivalent vertices as before.

\subsection{Second Differentials}
\begin{figure}[htb]
\begin{center}
\includegraphics[width=3in]{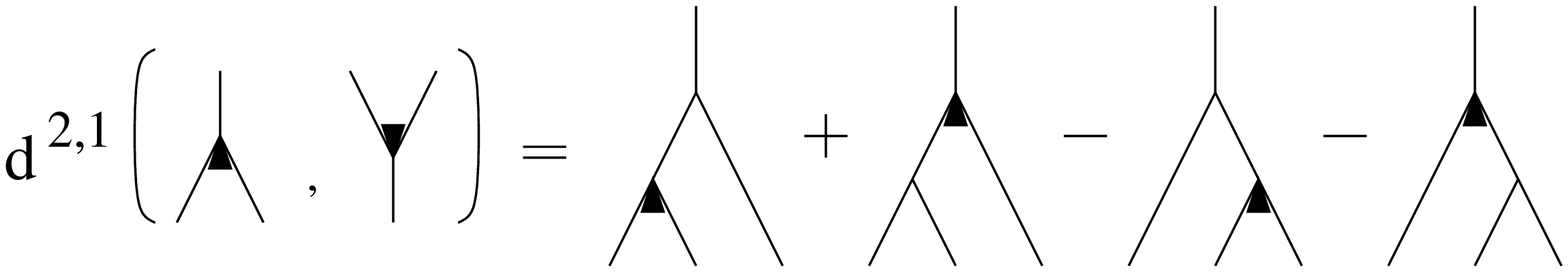}
\end{center}
\caption{A $2$-differential  $d^{2,1} $}
\label{d21}
\end{figure}

\begin{figure}[htb]
\begin{center}
\includegraphics[width=5in]{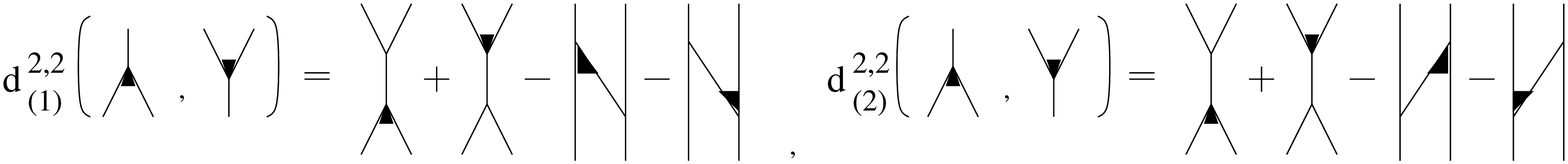}
\end{center}
\caption{$2$-differentials $d^{2,2}_{ (1) }$ and  $d^{2,2}_{ (2) }$}
\label{d22}
\end{figure}

An analogy with deformation theory when $\bar{\mu}=\mu + t \phi_1$
and $\bar{\Delta}=\Delta + t \phi_2$ gives the following.
\begin{definition}{\rm
Define the second  differentials by:
\begin{eqnarray*}
d^{2,1} (\phi_1 , \phi_2 ) = d^{2,1}(\phi_1) &=& \mu (\phi_1 \otimes |) + \phi_1 (\mu \otimes |)
- \mu (| \otimes \phi_1)
- \phi_1 (| \otimes \mu) , \\
d^{2,2}_{(1)} (\phi_1 , \phi_2) &=& \Delta \phi_1 + \phi_2 \mu
- (\phi_1 \otimes |)(| \otimes \Delta)
-
(\mu \otimes |)(| \otimes \phi_2), \\
d^{2,2}_{ (2)} (\phi_1 , \phi_2) &=& \Delta \phi_1 + \phi_2 \mu
- (| \otimes \phi_1)( \Delta \otimes |)
-
(| \otimes \mu)( \phi_2 \otimes |), \\
d^{2,3} (\phi_1 , \phi_2) = d^{2,3}(\phi_2) &=&  (\phi_2 \otimes |)\Delta
+  (\Delta \otimes |) \phi_2
- (| \otimes \phi_2) \Delta
- (| \otimes \Delta) \phi_2 .
\end{eqnarray*}
} \end{definition} Diagrams for $2$-cochain and $2$-differentials
are  depicted in Fig.~\ref{d21} for $d^{2,1} $ and  Fig.~\ref{d22} for $d^{2,2}_{ (1) }$ and  $d^{2,2}_{ (2) }$, respectively,
where $\phi_1 \in C^{2,1}_f(A;A)$ and $\phi_2\in C^{2,2}_f(A;A)$ 
are represented by black triangles on two different trivalent vertices.
The diagrams for $d^{2,3} $ are upside-down pictures of
Fig.~\ref{d21}.

Then define $D_2: C^2_f (A;A)  \rightarrow C^2_{f} (A;A)$ by
$D_2=D_2^{(i)}=d^{2,1}+d^{2,2}_{(i)}+d^{2,3}$ for either
$i=1$ or $2$.
Recall that for higher dimensions, differentials depend on this
choice of $i=1$ or $2$ at this exact dimension $2$ for $d^{2,2}_{
(1) }$ or $d^{2,2}_{ (2) }$, due to compatibility.
To avoid duplication in exposition, we choose,
once and for all, $i=1$. The case for $i=2$ will be clear, as all
the maps corresponding to $ (\mu \otimes |)(| \otimes \Delta)$ in
the case $i=1$ are simply replaced by those corresponding to $(|
\otimes \mu) (\Delta \otimes |)$ in the case $i=2$. Diagrammatically they are mirror images.
In the case $i=1$, the map looks like the letter ``N''
and its mirror in the case $i=2$.

\begin{figure}[htb]
\begin{center}
\includegraphics[width=4in]{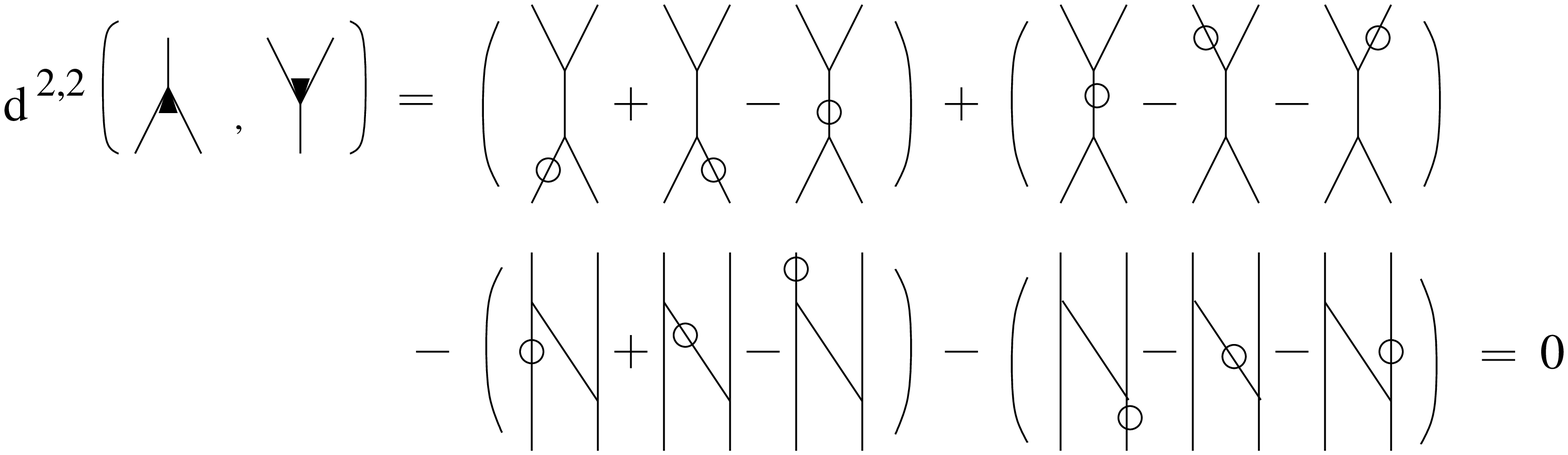}
\end{center}
\caption{The  $2$-cocycle condition for a $2$-coboundary:
$d^{2,2}_{ (1) }( d^{1,1}, d^{1,2} )=0$}
\label{dd22}
\end{figure}

\begin{theorem}\label{dd2thm}
$D_2 D_1=0$.
\end{theorem}
{\it Proof.\/}
This follows from direct calculations and also can be seen from diagrams as depicted in
Fig.~\ref{dd22}.
$\Box$

{}From the definition of the primary obstruction in Remark~\ref{deformrem},
the definition of the second differential implies the following.

\begin{lemma} \label{deform2lem}
For a deformation $(\bar{\mu}+\mu_{n+1}, \bar{\Delta}+\Delta_{n+1})$
of the multiplication and comultiplication of
a Frobenius algebra $A$ defined in Section~\ref{deformsec}
and the primary obstruction $(\xi_1, \xi_2, \xi_2', \xi_3)$ defined in
Remark~\ref{deformrem},
the relations $D_2^{(1)}(\mu_{n+1}, \Delta_{n+1})=(\xi_1, \xi_2, \xi_3)$
and  $D_2^{(2)}(\mu_{n+1}, \Delta_{n+1})=(\xi_1, \xi_2', \xi_3)$ hold.

In particular, the primary obstructions vanish if and only if
$(\mu_{n+1}, \Delta_{n+1})$ is a $2$-cocycle.
\end{lemma}

\subsection{Third differentials}

For simplicity we continue to use the case $d^{2,2}_{(1)}$ corresponding to
$ (\mu \otimes |)(| \otimes \Delta)$
for higher dimensions, and omit subscripts
$(1)$ from all the differentials.

\begin{definition}{\rm
Define the third  differentials by:
\begin{eqnarray*}
d^{3,1} (\xi_1 , \xi_2, \xi_3 ) &=& \mu (\xi_1 \otimes |) + \xi_1 (| \otimes \mu \otimes |)
+ \mu (| \otimes \xi_1)
- \xi_1 ( \mu \otimes | ) - \xi_1 ( |  \otimes  \mu)  , \\
d^{3,2}(\xi_1 , \xi_2, \xi_3) &=&
\Delta \xi_1 + \xi_2 (| \otimes \mu ) + ( \mu \otimes |)(| \otimes \xi_2 )
-\xi_2 (\mu  \otimes | ) - (\xi_1  \otimes |)( |^{\otimes 2} \otimes \Delta ),  \\
d^{3,3}(\xi_1 , \xi_2, \xi_3) &=&
\xi_3 \mu + (| \otimes \Delta ) \xi_2 + (| \otimes \xi_2 )( \Delta \otimes |)
- (\Delta  \otimes | )\xi_2 -( |^{\otimes 2} \otimes \mu ) (\xi_3  \otimes |),  \\
d^{3,4} (\xi_1 , \xi_2, \xi_3 ) &=&(\xi_3 \otimes |) \Delta +  (| \otimes \Delta \otimes |)\xi_3
+ (| \otimes \xi_3) \Delta
- \xi_1 ( \mu \otimes | ) -  ( |  \otimes  \Delta)\xi_3  .
\end{eqnarray*}
} \end{definition}
and recall that $$D_3=d^{3, 1} +d^{3, 2} +d^{3,
3} +d^{3, 4}:
 C^3_f (A;A)  \rightarrow C^4_f (A;A).$$

\begin{figure}[htb]
\begin{center}
\includegraphics[width=2.5in]{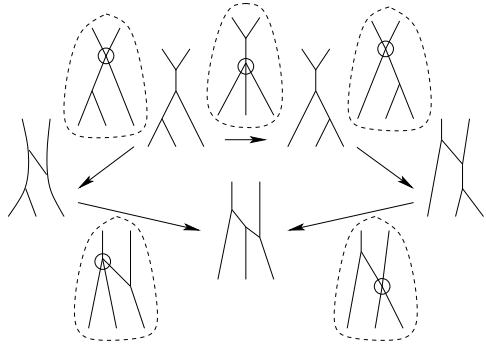}
\end{center}
\caption{Deriving a  $3$-differential $d^{3,2}$}
\label{d32circle}
\end{figure}

\begin{figure}[htb]
\begin{center}
\includegraphics[width=3.5in]{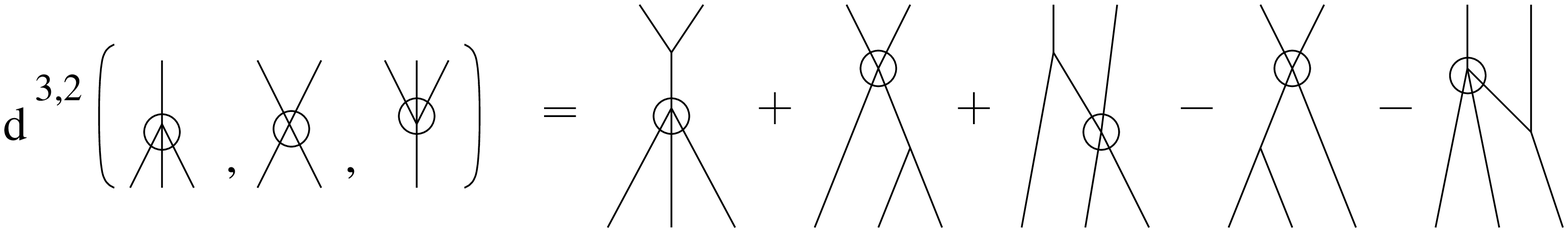}
\end{center}
\caption{A $3$-differential $d^{3,2}$}
\label{d32}
\end{figure}

The differentials $d^{3,1}$ and $d^{3,4}$ are the Hochschild differentials
that are derived from the pentagon conditions for associativity and coassociativity, respectively.
The differentials $d^{3,2}$ and $d^{3,3}$ are dual in the sense that once one is described via graphs, then the other is a mirror reflection in a horizontal line. We turn to describe the third differential $d^{3,2}$.  Recall that the $2$-cocycle conditions for a pair of maps $(\phi_1,\phi_2)$ are obtained by replacing incidences of multiplication and comultiplication by $\phi_1$ and $\phi_2$ in the associativity, compatibility, and coassociativity conditions. The equalities are depicted as equalities of graphs that contain trivalent vertices. The equalities themselves, then, can be depicted as 4-valent vertices. There are three such vertices: a $(3,1)$ vertex corresponds to associativity (3 up, 1 down), a $(2,2)$ vertex corresponds to either one of the compatibility conditions, and a $(1,3)$ vertex corresponds to coassociativity.
At the top left of  Fig.~\ref{d32circle} a graph is depicted that represents the composition $(\Delta)(\mu)(\mu \otimes |)$.  There are two ways of deforming this graph using associativity or the $N$-compatibility condition. These transformations are depicted in the
figure  as
an
encircled $(3,1)$ vertex on the right arrow or
an encircled $(2,2)$
vertex on the left arrow.
Such vertices represent $\xi_1$ or $\xi_2$, respectively.
Continuing around the diagram by choosing a new point
at which an identity could be applied,
we obtain a pentagon condition which
can be derived from the pentagon condition
on parentheses structures
that gives the Biedenharn-Elliot identity (see \cite{CFS}, for example).
Now read clockwise in a cycle around the diagram to write down the $3$-coboundary of the pair $(\xi_1,\xi_2)$; motion against the direction of the arrows results in negative coefficients.

\begin{figure}[htb]
\begin{center}
\includegraphics[width=4in]{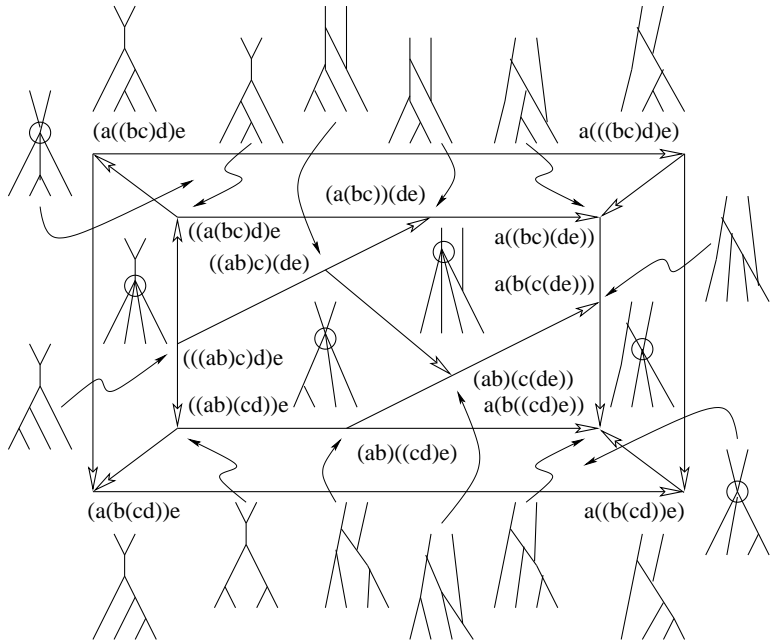}
\end{center}
\caption{Diagrams for $d^{4,2}$}
\label{assoc42}
\end{figure}

The pentagon itself can be labeled by a $(3,2)$-vertex that represents a $4$-cochain $\zeta_2 \in C^{3,2}_f(A;A).$ The corresponding $4$-coboundary  is a linear combination of the cells in the Stasheff polyhedron depicted in Fig.~\ref{assoc42} in which the cell corresponding to $\zeta_2$ appears as the bottom pentagon, which is drawn as a trapezoid in the figure.

\begin{figure}[htb]
\begin{center}
\includegraphics[width=4.5in]{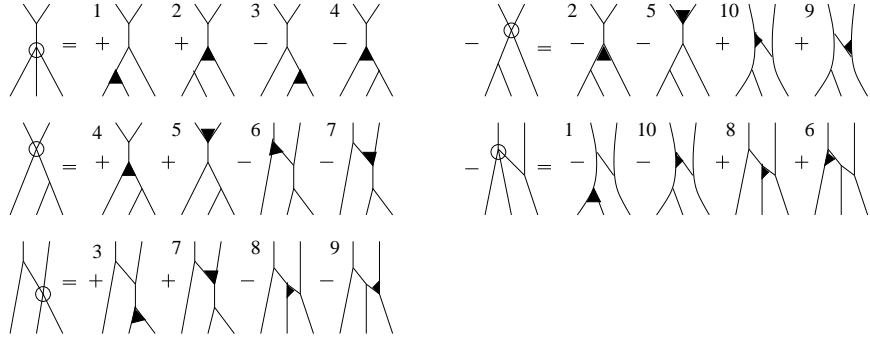}
\end{center}
\caption{$D_3 D_2=0$}
\label{dd32}
\end{figure}

\begin{theorem}\label{dd3thm}
$D_3 D_2=0$.
\end{theorem}
{\it Proof.\/} Since $d^{3,1}$ and $d^{3,4}$ are the same as the
Hochschild differentials of bialgebras, and $d^{3,2}$ and $d^{3,3}$
are dual to each other, we only consider the case of
$d^{3,2}[D_2(\phi_1, \phi_2)]=0$.
This again follows from direct calculations,
and can also be seen from diagrams as in Fig.~\ref{dd32}. In the
figure, each term of $d^{3,2}$ is replaced by the terms of
$D_2(\phi_1, \phi_2)$, and the canceling terms are indicated by
matching integers at the top-left corners of the diagrams. $\Box$

{}From point of view of deformation theory, continuing from
Lemma~\ref{deform2lem}, we obtain the following.

\begin{lemma} \label{deform3lem}
The primary obstruction $(\xi_1, \xi_2, \xi_3)$ and $(\xi_1, \xi_2',
\xi_3)$ to a deformation of the multiplication and comultiplication
of a Frobenius algebra $A$ defined in Remark~\ref{deformrem} give a
$3$-cocycle: $D_3(\xi_1, \xi_2, \xi_3)=0$, and
$D_3(\xi_1, \xi_2',
\xi_3)=0$.
\end{lemma}

\subsection{Fourth  differentials}

Recall that two third differentials, $d^{3,1}$ and $d^{3,4},$ are
the same as the Hochschild
differentials of bialgebras, and $d^{3,2}$ and $d^{3,3}$ are dual
to each other.
In fact, even the latter two have the same
diagrammatic aspect of being a pentagon. Indeed, the graphs
involved have two edges pointing upward instead of a single edge.
Since the same phenomena continues to the fourth differentials, we
only describe
the formula for $d^{4,2}$. As before, $d^{4,1}$ and $d^{4,5}$ are
Hochschild differentials, and $d^{4,4}$ is dual to $d^{4,2}$, and
$d^{4,3}$ is a self-dual symmetric version of $d^{4,2}$. Thus an
explicit formula, together with diagrammatic explanations for
$d^{4,2},$ would suffice.

\begin{definition} {\rm
Define $d^{4,2}$  as follows:
\begin{eqnarray*}
d^{4,2}(\zeta_1, \zeta_2, \zeta_3, \zeta_4)
&=& \Delta \zeta_1
+ \zeta_2 ( \id \otimes \id \otimes \mu)
+ (\mu \otimes \id)( \id \otimes \zeta_2) \\
& &  - (\zeta_1 \otimes \id)( \id \otimes \id \otimes \id \otimes \Delta)
- \zeta_2 ( \mu  \otimes \id \otimes \id )
- \zeta_2 ( \id \otimes \mu \otimes \id) .
\end{eqnarray*}
}
\end{definition}

Figure~\ref{assoc42} depicts the Stasheff polytope (or associahedron) with vertices, edges, and faces having labels adapted to the current purposes: to formulate the fourth differential.
The polytope consists of six pentagons and three squares. In the current drawing, four of the pentagons are drawn as trapezoids with their fifth vertex appearing as an interior  point on each of the shorter of the parallel edges. The boundary of the planar figure represents one of the three squares in the polytope. The drawing, then, is a distortion of a central projection through this face. The other two squares are in the northwestern and southeastern
center of the figure while the other two pentagons are in the southwest/northeast corners of the center.
To each vertex, 
we associate a parenthesization of the five letters $a,b,c,d,e$ 
so that the resulting string represents a composition of binary products. Each edge corresponds to a single
change of parenthesis such as $(ab)c$ to $a(bc)$. Meanwhile, a parenthesis structure can also be represented
by a tree diagram. For example, $(ab)c$ and $a(bc)$ are represented by the first and third
terms drawn in the RHS  in Fig.~\ref{d21}, respectively, where for the moment you should ignore the black triangles.  The pentagons in Fig.~\ref{assoc42} are formed as a cycle of such regroupings. The edges, vertices, and faces
here
are indicated by tree diagrams which have two upward pointing branches  and four downward pointing roots. To form  correspondences between these tree and the groupings of $a,b,c,d,e$ bend the right branches down.

The pentagon in Fig.~\ref{d32circle} represents a $4$-cochain
$\zeta_2 \in C^{4,2}(A;A)$, which is also represented by a
$5$-valent circled vertex with three edges pointing down and two
up. When  a $4$-cochain $\zeta_2$ is identified with a pentagon,
it is regarded as a homotopy of a path that starts from the three
edges that are three consecutive arrows of the same direction,
and sweeps the pentagonal face to the two remaining edge arrows. In
Fig.~\ref{d32circle}, the three
starting edges and two terminal edges both start at the top left
corner of the pentagon and
end at the bottom center.
If a homotopy starts from two arrows instead, and
sweeps the pentagon and ends at the three arrows, then it is
regarded as the negative of the corresponding cochain, $-\zeta_2$.

In Fig.~\ref{assoc42}, there is a unique vertex $v_0$ from which all three arrows
point out, which is assigned the
parenthesized term $(((ab)c)d)e$. The unique vertex $v_1$ into which three arrows point is
assigned
$a(b((cd)e))$. There are  edge paths that follow the directions of arrows that go
from $v_0$ to $v_1$.
Such an edge path $\gamma$ is homotopic to itself through
pentagonal and square faces (including the ``outside square''
in Fig.~\ref{assoc42}) of the Stasheff polytope
sweeping the entire sphere of the polytope.
By formulating each pentagonal face in terms of the maps represented
by tree diagrams, we obtain the formula for $d^{4,2}$.

\begin{theorem}\label{dd4thm}
$D_4 D_3=0$.
\end{theorem}
{\it Sketch of proof.\/} This again follows from direct calculations and can
also be seen from diagrams  in Fig.~\ref{assoc42}. Here we explain
how we see cancelations on diagrams. In computing $D_4 D_3$, the
terms in expressions in $D_4$  have their $\zeta$-factors
systematically replaced by  expressions that appear in the edges
of the boundary of each pentagon, but each edge is labeled by an
operator involving an $\xi$ as indicated, for example, in
Fig.~\ref{d32circle}. Since each edge is the boundary of exactly
two regions in the associahedron, the terms cancel. Note that if
one of the bounded cells is a square, then there is a corresponding
identity among those four terms --- commutativity of distant
tensor operators --- that causes the boundaries to cancel. The
remaining details are left to the reader. $\Box$

\begin{remark}{\rm While we have not followed through all the details of such a construction, it seems reasonable to parametrize all the higher differentials in terms of the cells of the higher-dimensional Stasheff polytopes. The coboundaries are parametrized by the boundaries of these cells, and that the square of the differential is trivial will follow from the codimension $1$ boundaries appearing on exactly two faces with opposite orientations. In this way, the cohomology of Frobenius algebras should be defined in all dimensions.}\end{remark}

\subsection{Cohomology Groups}

For convenience define
$C^0_f(A;A)=0$ and  $D_0=0: C^0_f(A;A)\rightarrow C^1_f(A;A)$.
Then Theorems~\ref{dd2thm} and \ref{dd3thm} are summarized as:
\begin{theorem}
${\cal C}=(C^n, D_n)_{n=0,1,2,3, 4}$ is a chain complex.
\end{theorem}

This enables us to define:
\begin{definition}{\rm
The {\it Frobenius $n$-coboundary, cocycle}, and {\it cohomology groups} are defined by:
\begin{eqnarray*}
B^{n}_f(A;A) &=& {\rm Image}(D_{n-1}) , \\
Z^{n}_f(A;A) &=& {\rm Ker}(D_n) , \\
H^{n}_f(A;A) &=& Z^{n}_f(A;A) / B^{n}_f(A;A)
\end{eqnarray*}
for $n=1,2,3, 4$.
} \end{definition}

The
lemma below
follows  from the definitions.

\begin{lemma} \label{d1lem}
Let $A$ be a Frobenius algebra
with $1=1_A=\eta{1_k}$.
Then:

\noindent
{\rm (i) }
$d^{1,1}(h)(1 \otimes 1) = h(1)$,  and

\noindent {\rm (ii) } $d^{1,1}(h)(1 \otimes x)=d^{1,1}(h)(x \otimes
1)=h(1)x,\;$
for any $x \in A$.

\noindent {\rm (iii) } If $\gamma(1)=1 \otimes x + x \otimes 1 $ for
some $x \in A$, and $h \in Z^1_f(A;A)$, then $h(x)=\alpha \cdot 1, $
for some constant $\alpha$ such that $2 \alpha =0 \in k$. In
particular, $h(x)=0$ if $\ch(k)\neq 2$.
\end{lemma}
{\it Proof.\/} One computes $d^{1,1}(h)(1 \otimes 1) = h(1)$,
and
(ii) follows from direct calculations. For (iii), using
$d^{1,2}(h)(1)=0$, we obtain $h(x)\otimes 1 + 1 \otimes h(x)=0$,
which implies the statement.  $\Box$

\begin{lemma} \label{d2lem}
If a Frobenius algebra $A$  is commutative and $d^{2,1}(\phi_1)=0$,
then for any $x \in A$, the following hold:
\begin{eqnarray*}
& & x \phi_1(1\otimes 1)=\phi_1(1\otimes x)=\phi_1(x\otimes 1) , \\
& & \phi_1(x^2\otimes x)=\phi_1(x\otimes x^2),  \quad
\phi_1(1\otimes x^2)=x \phi_1(1 \otimes x) .
\end{eqnarray*}
\end{lemma}
{\it Proof.\/} One computes $d^{2,1}(\phi_1)(a \otimes b \otimes c) $
for
$(a,b,c)=(1,1,x)$ and $ (x,1,1)$,
for the first set of equations, and
$(x, x,x)$, $(1,x,x)$ for the second set,
respectively. The other choices using two elements
$\{ 1, x \}$ do not give additional conditions. $\Box$

\subsection{Examples}

In this subsection we choose $d^{2,2}_{(1)}$ for the chain complex to compute.
Throughout this section, the symbols $\gamma^{ab}_c$ and $\lambda^{a}_{bc}$ indicate the structure constants of comultiplication and multiplication for the algebras in question.

\begin{example}\label{cplxcoh}{\rm
For the example of complex numbers in Example~\ref{cplxex},
we have
$$H^1_f (\C;\C)  =  0,  \quad
Z^2_f (\C ; \C) = \R^6, \quad
H^2_f (\C ; \C) =  \R^4 .
$$

\noindent
{\it Proof.\/}
By Lemma~\ref{d1lem} (i), for $h \in Z^1_f (\C; \C)$,
we have $h(1)=0$, and from
$d^{1,1}(h)(i \otimes i)=0$, we obtain $h(i)=0$.
Thus we obtain $Z^1_f(\C; \C)=0=H^1_f (\C;\C) $.
This also implies that  $B^2_f(\C; \C)\cong \R^2$.

By Lemma~\ref{d2lem},
we have
$$ i \phi_1(1\otimes 1)=\phi_1(1\otimes i)=\phi_1(i\otimes 1) $$
from $d^{2,1}(\phi_1)=0$.
Hence we write
$$\phi_1(1\otimes i)=\phi_1(i\otimes 1)
= \lambda 1 + \lambda' i $$ for some $\lambda, \lambda' \in \R$, and
then $ \phi_1(1\otimes 1)= \lambda' 1 - \lambda i $.
By setting
$\phi_2(a)=\sum_{b,c} \gamma_a^{b,c} (b \otimes c)$ for basis
elements $a, b, c \in \{ 1, i\}$, $d^{2,3}(\phi_2)=0$ implies
\begin{eqnarray}
\gamma_1^{1,1} = \gamma_i^{1,i}=\gamma_i^{i,1}, & & \label{Ccocy2eq1} \\
- \gamma_i^{1,1} = \gamma_1^{1,i}=\gamma_1^{i,1}, & & \label{Ccocy2eq2}
\end{eqnarray}
leaving free variables $\gamma_1^{1,1}$,  $\gamma_1^{i,i}$,
$\gamma_i^{1,1}$ and $\gamma_i^{i,i}$. The equations
$d^{2,2}_{(1)}(\phi_1, \phi_2)(1 \otimes 1)=0$ and
$d^{2,2}_{(1)}(\phi_1, \phi_2)(1 \otimes i)=0$ do not give
additional conditions.  Assuming that $\phi_1(i \otimes
i)=\lambda_{i,i}^1\;1 + \lambda_{i,i}^i\;i,$
the other two conditions of $d^{2,2}_{(1)}(\phi_1, \phi_2)=0$
give
$$ - \lambda + \lambda_{i,i}^i  + \gamma_i^{i,i} - \gamma_1^{1,i}= 0,
\quad
\lambda' + \lambda_{i,i}^1  + \gamma_1^{i,i} +\gamma_i^{1,i}=0 ,
$$
that 
can be rewritten with Eqns.~(\ref{Ccocy2eq1}) and (\ref{Ccocy2eq2})
as
\begin{eqnarray}
- \lambda + \lambda_{i,i}^i  + \gamma_i^{i,i} + \gamma_i^{1,1}&=& 0,
\label{Ccocy2eq3} \\
\lambda' + \lambda_{i,i}^1  + \gamma_1^{i,i} +\gamma_1^{1,1}&=&0 .
\label{Ccocy2eq4}
\end{eqnarray}
Since  dim$(Z^2_f(A;A))$  is equal to the number of variables
$( \lambda, \lambda',   \lambda_{i,i}^1,  \lambda_{i,i}^i,
\gamma_1^{1,1}, \gamma_1^{i,i}, \gamma_i^{1,1}, \gamma_i^{i,i}   )$
minus
the number of equations (Eqns.~(\ref{Ccocy2eq3}) and (\ref{Ccocy2eq4})),
we obtain dim$(Z^2_f(A;A))=6$.
Hence, together with $B^2_f(\C; \C)\cong \R^2$,
we obtain the stated results.

} \end{example}

\begin{example}\label{polycoh}{\rm
For $A=k[x]/(x^2)$  in Example~\ref{polyex}, we have
$$ H^1_f (A;A) = \left\{
\begin{array}{l} 0 \quad {\rm if } \quad
\ch(k) \neq 2  \\ k  \quad {\rm if }   \quad
\ch(k) = 2
\end{array} \right. ,
\quad
Z^2_f (A; A)
= k^6 ,
\quad
H^2_f (A;A)
=
\left\{
\begin{array}{l} k^4 \quad {\rm if }  \quad  \ch(k) \neq 2
\\ k^5 \quad {\rm if }  \quad  \ch(k) = 2
\end{array} \right.  
$$
{\it Proof.\/}
{}From the proof of  Lemma~\ref{d1lem}, the condition
$h \in Z^1_f(A;A)$ is equivalent to
$h(1)=0$, $h(x)=\alpha \cdot 1$ with $2 \alpha =0,$  and
the following additional conditions that were not used in the proof:
\begin{eqnarray*}
d^{1,1}(h)(x\otimes x) & = & h(x) \cdot x + x \cdot h(x) - h(0)
=2xh(x)=0, \\
d^{1,2}(h)(x)&=& h(x) \otimes x +  x \otimes h(x)
- \Delta (h(x))=0,
\end{eqnarray*}
both of which follow from the conditions already stated in the
lemma. Hence we obtain $H^1_f$ as stated. We also have
$B^2(A;A)\cong k^2$ if $\ch(k)\neq 2$ and $B^2(A;A)
\cong k$  if $\ch(k)=2$.

For  $2$-cocycles $\phi_1 \in C^{2,1}(A;A)$ and $\phi_2 \in C^{2,2}(A;A)$,
Lemma~\ref{d2lem} implies that there is $\lambda \in k$ such that
$$\phi_1(1 \otimes x)=\phi_1(x \otimes 1)= \lambda x$$ and $\phi_1(1 \otimes 1)=\lambda 1 + \lambda' x$ for another
$\lambda' \in k$. Direct calculations show also if $\phi_2 (a) =
\sum_{b, c} \gamma_a^{b, c} (b \otimes c)$ then $d^{2,3}(\phi_2)=0$ implies
$$ \gamma_x^{1,x}= \gamma_x^{x,1}=0,
\quad \gamma_1^{1,x}= \gamma_1^{x,1}=\gamma_x^{x,x}. $$ Now let
$\phi_1(x \otimes x)=\alpha 1 +\beta x$.  The equation
$d^{2,2}_{(1)}(\phi_1, \phi_2)(x \otimes 1)=0$ implies  $\phi_1(x
\otimes x)=\gamma_x^{1,1} 1-\gamma_1^{1,1} x$, (the evaluations at
other tensors don't give any extra conditions).  In summary we
obtain
$\phi_2(x)=\gamma_x^{1,1}(1 \otimes 1)+ \gamma_1^{1,x}(x \otimes x)$ and
$\phi_2(1)= \gamma_1^{1,1} (1 \otimes 1) +\gamma_1^{1,x} ( 1 \otimes x +x \otimes 1) + \gamma_1^{x,x}(x \otimes x )$,
in total  a six-dimensional solution set parametrized by $\lambda, \lambda', \gamma_1^{1,1}, \gamma_x^{1,1},\gamma_1^{1,x}$ and $\gamma_1^{x,x}$.
The result follows.
} \end{example}

\begin{example}\label{gpcoh}{\rm
For a group algebra $A=kG$ in Example~\ref{gpex},
we consider the case $G=\Z_2$. Then we have
$$
H^1_f (A;A) =
\left\{
\begin{array}{l} 0 \quad {\rm if }  \quad \ch(k) \neq  2  \\
           k  \quad {\rm if } \quad \ch(k) = 2
\end{array} \right. ,
\quad
Z^2_f (A; A)= k^6 ,
\quad
H^2_f (A;A)
=
\left\{
\begin{array}{l} k^4  \quad {\rm if } \quad \ch(k) \neq 2  \\
           k^5  \quad {\rm if } \quad \ch(k) = 2
\end{array} \right.  
$$
{\it Proof.\/}
Assuming $d^{1,1}(h)=0$ for
$h \in  C^1_f(A;A)$, Lemma~\ref{d1lem} implies that
$h(1)=0$. The condition  $d^{1,1}(h)(x \otimes x)=0$ implies
$2xh(x)=0$, which is equivalent to $2h(x)=0$.
The same condition follows from  $d^{1,2}(h)(1)=0$,
and the last condition  $d^{1,2}(h)(x)=0$ implies that
$h(x)=\alpha x $ for some $\alpha \in k$.
Thus we obtain $H^1$ as stated.

Lemma~\ref{d2lem} implies that
$\phi_1$ is given by
$$\phi_1(1\otimes x)=\phi_1(x\otimes 1)
= \lambda 1 + \lambda' x $$
for some $\lambda, \lambda' \in k$, and
$ \phi_1(1\otimes 1)= \lambda' 1 + \lambda x $.

{} From $d^{2,3}(\phi_2)=0 $
we obtain
$$ \gamma_1^{1,1}=\gamma_x^{1,x}=\gamma_x^{x,1},
\quad  \gamma_x^{1,1}=\gamma_1^{1,x}=\gamma_1^{x,1}.$$
Hence, we can write
\begin{eqnarray*}
\phi_2(1)
&=& q (1 \otimes 1) + r (1 \otimes x + x \otimes 1)+ \gamma_1^{x,x}(x \otimes x)
,\\
\phi_2(x)
&=& r(1 \otimes 1) + q (1 \otimes x + x \otimes 1) + \gamma_x^{x,x} (x \otimes x).
\end{eqnarray*}
Now let $\phi_1(x \otimes x)=\lambda^1_{x,x}\;1+ \lambda^x_{x,x}\;x$.  
The equation $d^{2,2}(\phi_1, \phi_2)=0 $
gives by evaluation at the four basis elements the following constraints
 \begin{eqnarray}
\gamma_1^{x,x}+ \lambda_{x,x}^1 &=& q + \lambda' ,  \label{z2eqn1} \\
 \gamma_x^{x,x}- \lambda_{x,x}^x &=& r -  \lambda .  \label{z2eqn2}
 \end{eqnarray}
Since  dim$(Z^2_f(A;A))$  is equal to the number of variables $(q,r,
\lambda, \lambda ', \gamma_1^{x,x},  \gamma_x^{x,x},
\lambda_{x,x}^1, \lambda_{x,x}^x$)
minus the number of equations
(the above two), we obtain dim$(Z^2_f(A;A))=6$. Thus we obtain the
result.

}\end{example}

\begin{remark}{\rm It is interesting that these Frobenius algebras all have the same $1$ and $2$-dimensional cohomology when char($k$)$\ne2$.
Nevertheless, the free variables are quite a bit different in each computation. We expect that higher-dimensional algebras are cohomologically distinct.}
\end{remark}

\section{Yang-Baxter solutions in Frobenius algebras and their cocycle deformations}\label{YBEsec}

\begin{figure}[htb]
\begin{center}
\includegraphics[width=3in]{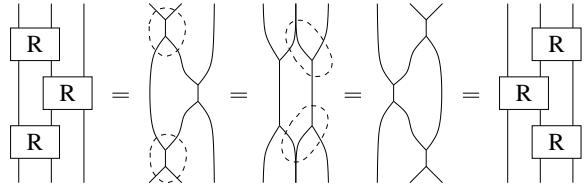}
\end{center}
\caption{A solution to YBE  in Frobenius algebras}
\label{YBE1}
\end{figure}

\begin{figure}[htb]
\begin{center}
\includegraphics[width=3.5in]{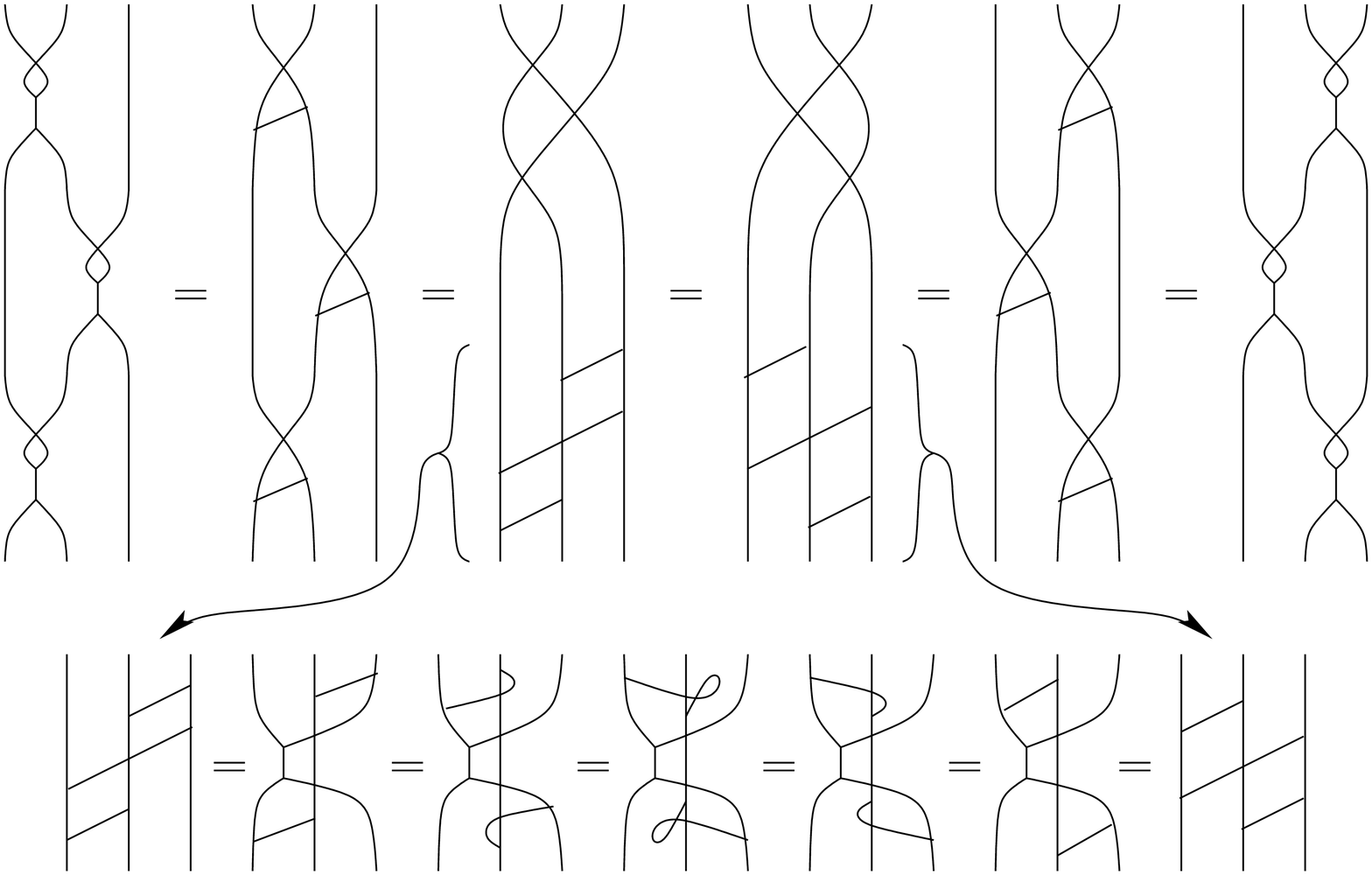}
\end{center}
\caption{A solution to YBE in symmetric Frobenius algebras}
\label{AmidaYBE}
\end{figure}

\begin{figure}[htb]
\begin{center}
\includegraphics[width=4in]{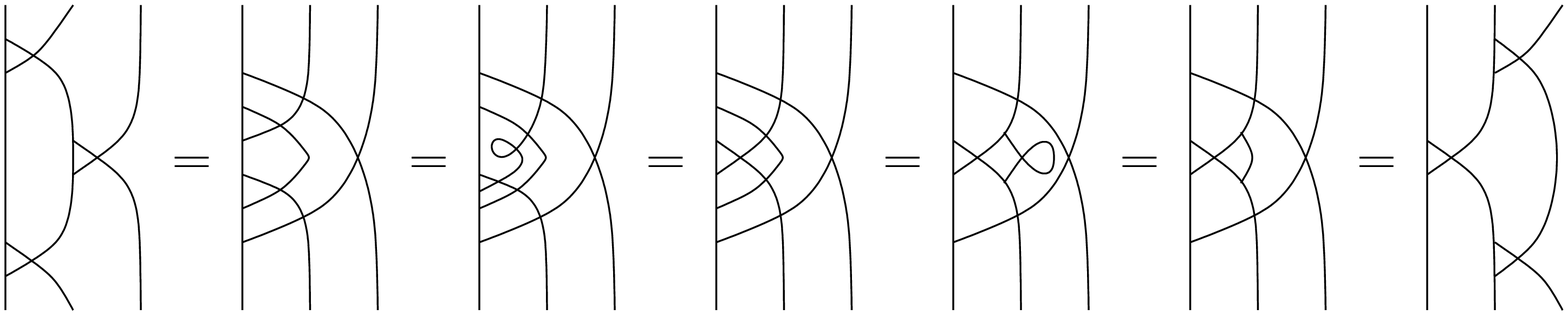}
\end{center}
\caption{Another solution to YBE in symmetric Frobenius algebras}
\label{KYBE}
\end{figure}

In this section, we construct YBE solutions ($R$-matrices)
from skein theoretic methods using maps in Frobenius algebras.
We start with the following examples, the first of which is due to
\cite{Stolin}, and diagrammatic proofs are depicted in
Figs.~\ref{YBE1}, \ref{AmidaYBE}, and  \ref{KYBE}, respectively.

The Yang-Baxter equation (YBE) is formulated as $(R \otimes
\id)(\id \otimes R)(R \otimes \id) =(\id \otimes R)(R \otimes
\id)(\id \otimes R)$ for $R \in {\rm Hom}(V \otimes V, V \otimes
V)$, see, for example, \cite{K&P}. It is often required that $R$
be invertible, but we do not always require the invertibility
unless explicitly mentioned. In particular, the solutions in
Lemma~\ref{ybe1lem} are not invertible. Invertible solutions
derived from these are discussed in Proposition~\ref{skeinprop}.

\begin{lemma}\label{ybe1lem}
For any Frobenius algebra $X$,
$R=\Delta \mu$ is a solution
to the Yang-Baxter equation.
For any symmetric Frobenius algebra $X$,
$R_1=\tau \Delta \mu$ and
$R_2 =(\mu \otimes \id )( \id \otimes \tau) (\Delta \otimes \id)$
are solutions
to the YBE as well.
\end{lemma}

\begin{example}{\rm
For a group ring $X=kG$ in Example~\ref{gpex},
the $R$-matrix  of type   $\Delta \mu$ is
computed by $R(x \otimes y) = (\Delta \mu) (x \otimes y)=
\sum_{zw=xy}z \otimes w$
for $x, y \in G$.
} \end{example}

\begin{figure}[htb]
\begin{center}
\includegraphics[width=3in]{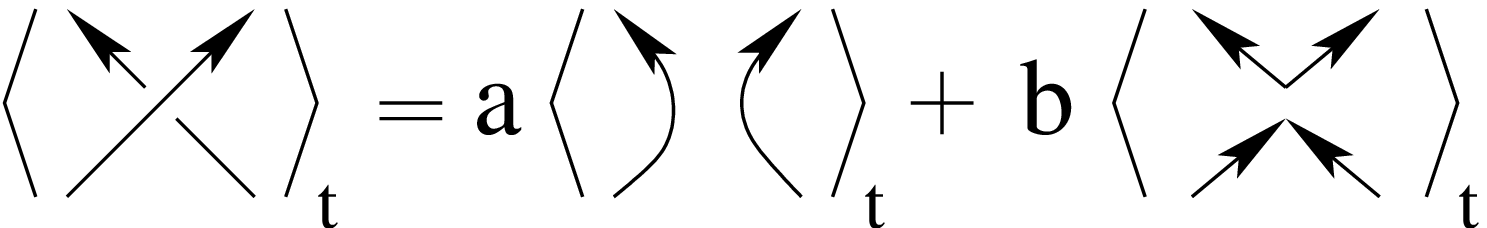}
\end{center}
\caption{A skein relation}
\label{skein}
\end{figure}

We define
two $R$-matrices $R$ and $R'$
by
skein relations as depicted in Fig.~\ref{skein},
which are written as follows. 
\begin{eqnarray}
R & = &  A \id^{\otimes 2} + B ( \gamma \beta ) + C (\Delta \mu) + T(\tau), \label{skein1}\\
R' &=&  A' \id^{\otimes 2} + B' (\gamma \beta ) + C' (\Delta \mu)+ T'(\tau). \label{skein2}
\end{eqnarray}
We call this the {\it Frobenius skein relation},
and below we derive conditions for it to satisfy the YBE.
We also suppose that the inverse of $R$ is
given by $R'$, and compute the conditions for this requirement.
Note that instead of the map $ \Delta \mu$,
the two other  solutions in Lemma~\ref{ybe1lem}
can be used to define similar skein relations, and
these possibilities might deserve further study.

Note that
$\beta \gamma (1)$ is an element of $k$ which we denote by $\delta_0$.
In the following, sometimes we make the assumption that
$\mu \Delta = \delta_1 {\rm \id}$ 
for some $\delta_1 \in k$
for computational simplicity.
(Recall  that in general it is true that
$\mu \Delta = \delta_h {\rm \id}$ 
holds but $\delta_h \in A$
is a central element.)
This holds for some of the examples
of Frobenius algebras, such as
group algebras.
Under this assumption one obtains
the following relation:
\begin{eqnarray}
\delta_0 =  \beta \gamma= \epsilon \mu \Delta \eta (1)
=\epsilon (\delta_1 1 ) = \delta_1 \epsilon (1). \label{delta0eq}
\end{eqnarray}

\begin{figure}[htb]
\begin{center}
\includegraphics[width=5in]{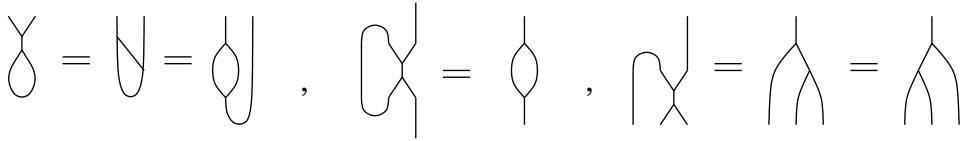}
\end{center}
\caption{A few formulas}
\label{bubble2}
\end{figure}

\begin{figure}[htb]
\begin{center}
\includegraphics[width=4.5in]{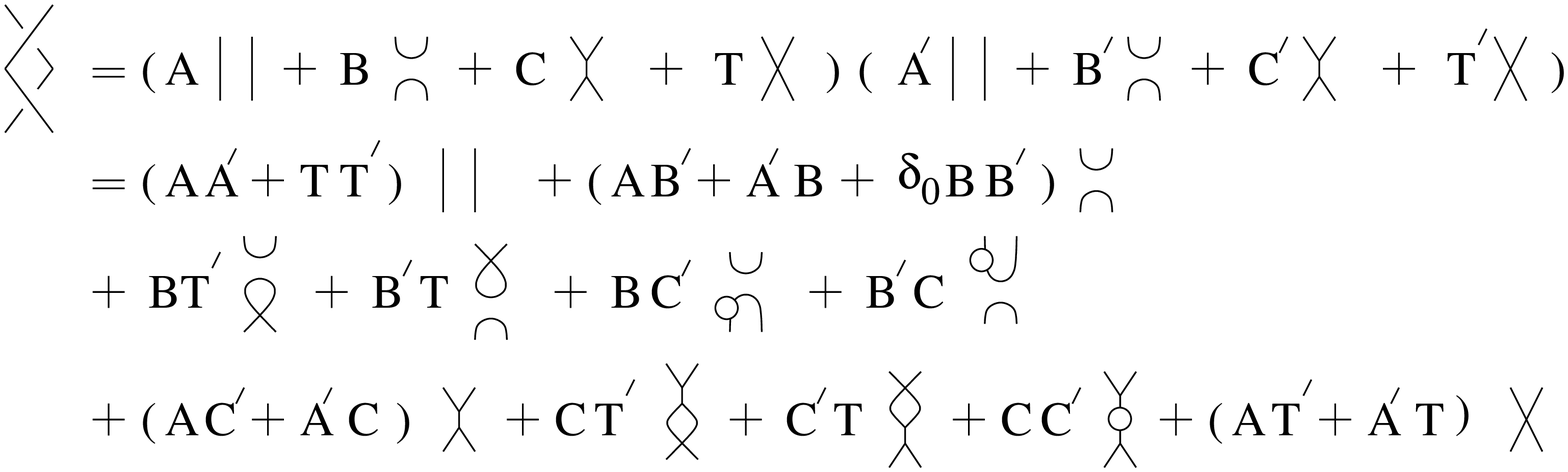}
\end{center}
\caption{The  skein relation of the inverse}
\label{skeininv}
\end{figure}

In  Fig.~\ref{bubble2}  a few direct calculations are depicted that will be used below.
Using 
this notation, one calculates the condition
for $R'$ to be the inverse of $R$.
A calculation is illustrated in Fig.~\ref{skeininv} for
the condition that $R'=R^{-1}$.

\begin{figure}[htb]
\begin{center}
\includegraphics[width=5.5in]{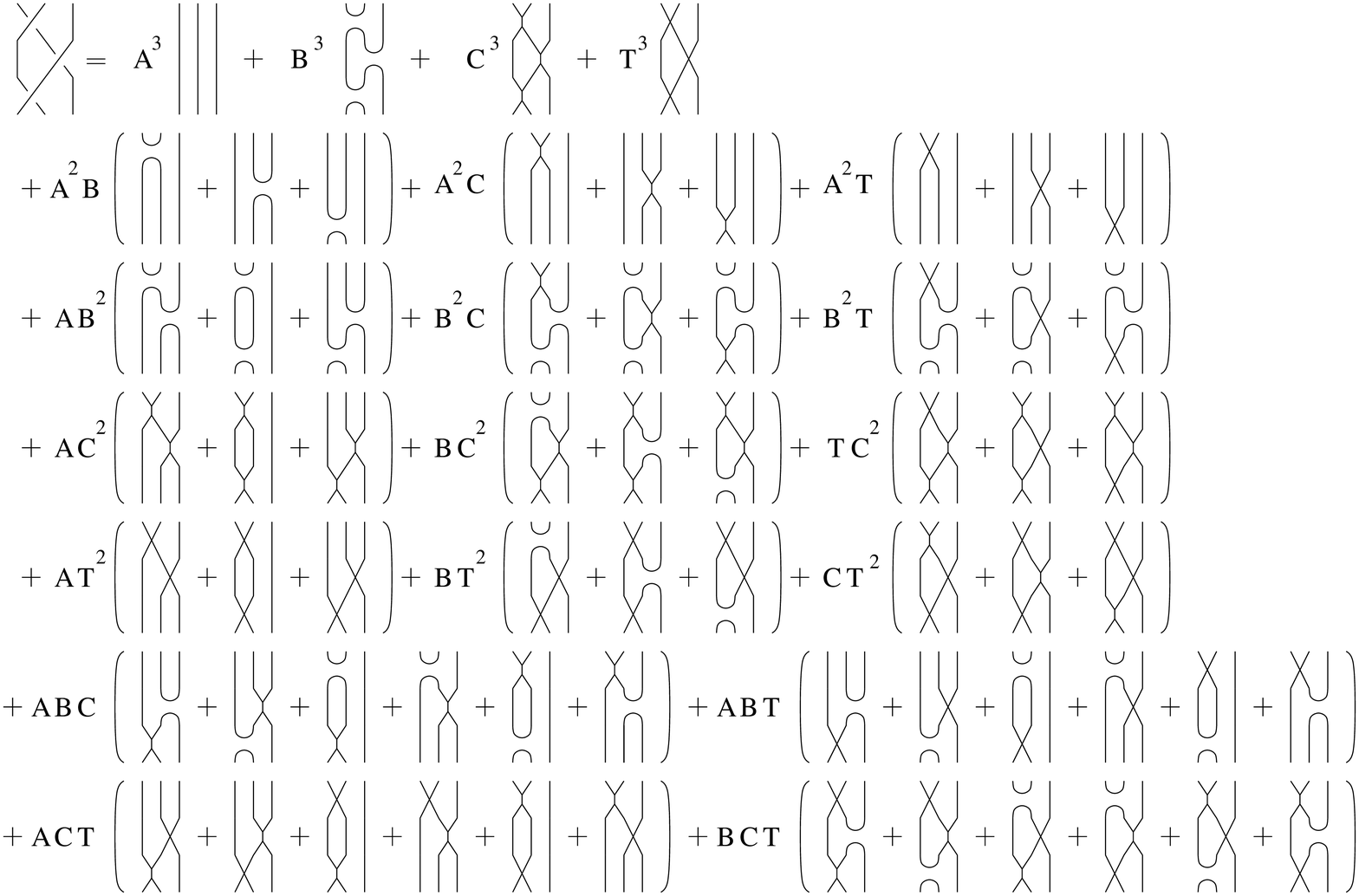}
\end{center}
\caption{YBE for a skein, LHS}
\label{skeinL}
\end{figure}

\begin{figure}[htb]
\begin{center}
\includegraphics[width=5.5in]{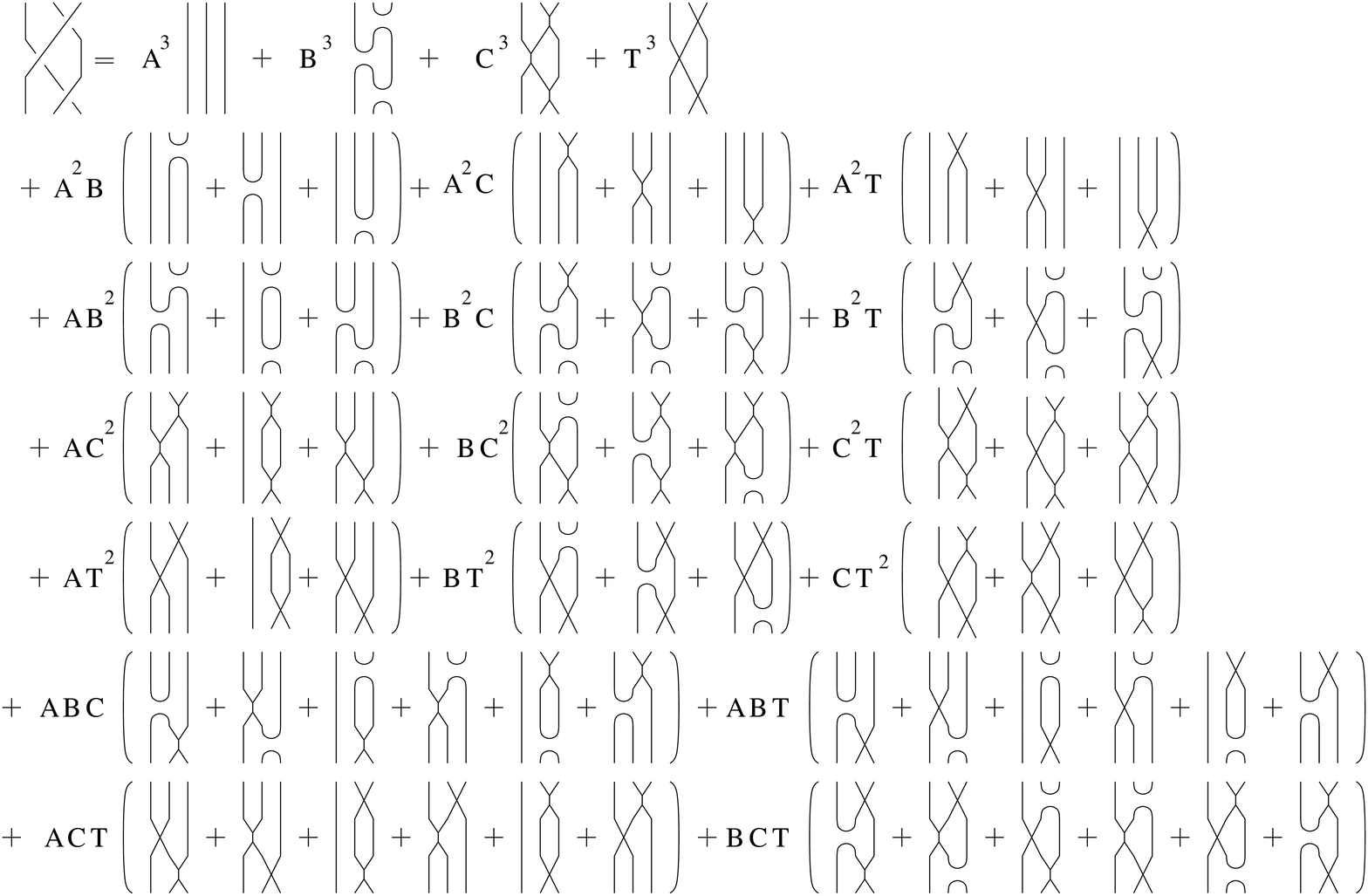}
\end{center}
\caption{YBE for a skein, RHS}
\label{skeinR}
\end{figure}

In Figs.~\ref{skeinL} and \ref{skeinR}
the YBE is formulated diagrammatically for the above defined skein relation.

\begin{proposition}\label{skeinprop}
Suppose
the Frobenius algebra $X$ over a field $k$
satisfies $\mu \Delta = \delta_1 {\id}$
for some $\delta_1 \in k$.
Then  the
$R$-matrix defined by the Frobenius skein relation~(\ref{skein1}),
with the inverse $R^{-1}=R'$ defined by the relation~(\ref{skein2}),
gives a solution to the YBE if
the following
hold:

\noindent
{\rm (i)}
$C=T=0$, $C'=T'=0$, 
$A^2 + B^2 + \delta_0 AB=0$,
$A'^2 + B'^2 + \delta_0 A'B'=0$, 
and $A B'+A'B+\delta_0 BB'=0$.

\noindent
{\rm (ii)}
$X$ is commutative,
$A=B=0$, $A'=B'=0$, $TT'=1$ and $CT'+C'T+ \delta_1 CC'=0$.
\end{proposition}
{\it Proof.\/}
We obtain the conditions
$A^2 T=0$ and $B^2 T=0$
by comparing the coefficients of the following terms for both sides of the equation
respectively (see Figs.~\ref{skeinL} and \ref{skeinR}):
\begin{eqnarray*}
& &\id \otimes \tau , \ \tau \otimes \id  , \\
& &(\id \otimes \tau)(\gamma \otimes \id)(\beta  \otimes \id),
\ (\tau \otimes \id)( \id \otimes \gamma)( \id \otimes \beta) ,
\ (\beta  \otimes \id)(\gamma \otimes \id)(\id \otimes \tau),
\ ( \id \otimes \beta)( \id \otimes \gamma)(\id \otimes \tau) .
\end{eqnarray*}
Assuming that variables take values in the field $k$, we obtain
either $T=0$ or $A=B=0$.

Assume $T=0$, and
also assume that $R^{-1}$ satisfies the similar condition,  $T'=0$.
For $R$ we obtain
\begin{eqnarray}
B( A^2 + B^2 + \delta_0 AB + \delta_1 BC + 2 \delta_1 A C  ) &=&0,  \label{skeineq1} \\
AC( A +  \delta_1 C ) &=&0,  \label{skeineq2}  \\
BC( B + \delta_1 C  ) &=&0,  \label{skeineq3}
\end{eqnarray}
that are derived  from comparing the coefficients  of the
following maps:
\begin{eqnarray*}
& &\id \otimes \gamma \beta, \   \gamma \beta \otimes \id , \\
& &(\id \otimes \Delta \mu) , \ (\Delta \mu \otimes \id) , \\
& &(\id \otimes \gamma)\mu (\mu \otimes \id) ,
\ ( \gamma \otimes \id )\mu (\mu \otimes \id) ,
\ (\Delta \otimes \id)\Delta (\id \otimes \beta),
\ (\id  \otimes \Delta )\Delta (\beta \otimes \id ) .
\end{eqnarray*}
It is, then,  easy to see from these equations and from
Figs.~\ref{skeinL} and \ref{skeinR}
that the conditions
$C=T=0$ and   $A^2 + B^2 + \delta_0 AB=0$
give solution $R$ to the YBE. Similarly from the condition
depicted in Fig.~\ref{skeininv},
it follows that $R'$ gives the inverse of $R$ with the conditions
$C'=T'=0$,   $A'^2 + B'^2 + \delta_0 A'B'=0$
and $A B'+A'B+\delta_0 BB'=0$.

In this paragraph we observe
that the condition $C=0$, in fact, follows from the
required conditions assuming that variables take values in a field $k$.
Suppose $C\neq 0$.
Then Eqns. (\ref{skeineq2}) and (\ref{skeineq3})
require $A=B=-\delta_1 C $ and $A'=B'=-\delta_1 C' $.
Then Eqn. (\ref{skeineq1}) implies $\delta_0=1$, and
from Eqn. (\ref{delta0eq}) ($\delta_0=\delta_1 \epsilon(1)$),
we have $\delta_1=\epsilon(1)^{-1}$.
The inverse formula (Fig.~\ref{skeininv})
in this case ($C=T=0$) requires
$A B'+A'B+\delta_0 BB'=0$, which reduces to $\delta_1CC'=0$,
a contradiction.
Hence we
have $C=0$.

Next we consider the case
$A=B=0$, and assume that $R^{-1}$ is in the same form, $A'=B'=0$.
Then the commutativity implies LHS$=$RHS in Figs.~\ref{skeinL} and \ref{skeinR},
using formulas depicted in Fig.~\ref{bubble2}.
The inverse condition in
Fig.~\ref{skeininv} and the definition $R^{-1}=C'(\Delta \mu )+ T' \tau$
imply $TT'= 1$ and
$CT'+C'T+\delta_1 CC'=0$.
This leads to Case (ii).
$\Box$

\bigskip

Using $2$-cocycles,
we construct new $R$-matrices from old by deformation
as follows.
Let $X$ be a Frobenius algebra over $k$.
Suppose $R$ is defined by the Frobenius skein relation
that satisfies the conditions in Proposition~\ref{skeinprop},
so that $R$ is a solution to the YBE on $X$.

Let $\hat{X}=A \otimes k[[t]] / (t^2)$.
Then $\hat{X}$ is regarded as $ (k[t]/ (t^2) )$-module.
Extend the maps $\mu$ and $\Delta$ to $\hat{X}$.
{}From the
deformation interpretation of $2$-cocycles in Section~\ref{deformsec},
we have the following.

\begin{theorem}\label{deformRthm}
Let $X$ be a commutative, and therefore cocommutative, Frobenius algebra.
Suppose $\phi_i \in C^{2}(X;X)$, $i=1,2$,
are Frobenius $2$-cochains satisfying all the
$2$-cocycle conditions
$d^{2,1} = d^{2,2}_{(1)} =d^{2,2}_{(2)}=d^{2,3} =0$.
Define $R_{\phi_1, \phi_2} : \hat{X} \otimes  \hat{X} \rightarrow  \hat{X} \otimes \hat{ X} $
by
$$R_{\phi_1, \phi_2}=
C ( (\Delta+ t \phi_2)( \mu + t \phi_1)) + T(\tau). $$
 Then
$R_{\phi_1, \phi_2}$ is a solution to the YBE if the following conditions are satisfied:
$(\mu+t \phi_1)(\Delta + t \phi_2)=\delta_1 \id$ on $ \hat{X}$
for some $\delta_1 \in   k[t]/ (t^2)$, 
$\phi_1 \tau=\phi_1$, and $\tau \phi_2 = \phi_2$.
\end{theorem}
{\it Proof.\/}  This is a repetition of the proof of
Proposition~\ref{skeinprop} 
Case (ii), using the deformation interpretations
of $2$-cocycle conditions. The associativity, coassociativity, and
Frobenius compatibility conditions for $\mu + t \phi_1$ and $\Delta +
t \phi_2$ follow from the $2$-cocycle conditions in
Lemma~\ref{deformlem}. The conditions
$\phi_1 \tau=\phi_1$ and $\tau \phi_2 = \phi_2$
correspond to the commutativity. $\Box$

\begin{example} \label{cplxRex}
{\rm
For $X=\C$ in Example~\ref{cplxex},
the general solutions for the  $2$-cocycles $\phi_1$ and $\phi_2$
with $d^{2,1} = d^{2,2}_{(1)} =d^{2,2}_{(2)}=d^{2,3} =0$ found in
Example~\ref{cplxcoh} satisfy the condition (ii)
in Theorem~\ref{deformRthm}: $\phi_1 \tau=\phi_1$, $\tau \phi_2 = \phi_2$.
We check the condition
$(\mu+t \phi_1)(\Delta + t \phi_2)=\delta_1 \id$ on $ \hat{X}$
for some $\delta_1 \in   k[t]/ (t^2)$.
One computes:
\begin{eqnarray*}
(\mu+t \phi_1)(\Delta + t \phi_2)(1)&=&
2+ t\ [  \ (\gamma_1^{1,1} - \gamma_1^{i,i} + \lambda' -\lambda_{i,i}^1)
+ ( \gamma_1^{1,i} + \gamma_1^{i,1} - \lambda -\lambda_{i,i}^i) i \ ] , \\
(\mu+t \phi_1)(\Delta + t \phi_2)(i)&=&
2i+ t\ [  \ (\gamma_i^{1,1} - \gamma_i^{i,i} +2 \lambda )
+ ( \gamma_i^{1,i} + \gamma_i^{i,1} + 2 \lambda' ) i \ ] .
\end{eqnarray*}
Thus the general $2$-cocycles  satisfy
$(\mu+t \phi_1)(\Delta + t \phi_2)=\delta_1 \id$
if and only if the above two values are
multiples of $1$ and $i$, respectively, by the same element $\delta_1 \in  k[t]/ (t^2)$.
This condition is written as
\begin{eqnarray}
\gamma_1^{1,i} + \gamma_1^{i,1} - \lambda -\lambda_{i,i}^i &=& 0 ,\label{CR1}\\
\gamma_i^{1,1} - \gamma_i^{i,i} +2 \lambda  &=& 0 ,\label{CR2} \\
 \gamma_1^{1,1} - \gamma_1^{i,i} + \lambda' -\lambda_{i,i}^1
&=& \gamma_i^{1,i} + \gamma_i^{i,1} + 2 \lambda' . \label{CR3}
\end{eqnarray}
For Eqn.~(\ref{CR3}), Eqns.~(\ref{Ccocy2eq1}) and (\ref{Ccocy2eq4})
imply
$$ \gamma_1^{1,1} - \gamma_1^{i,i} + \lambda' -\lambda_{i,i}^1=
-2 (\gamma_1^{i,i} +  \lambda_{i,i}^1)
= \gamma_i^{1,i} + \gamma_i^{i,1} + 2 \lambda' , $$
so that Eqn.~(\ref{CR3}) is redundant and we obtain
$\delta_1=2[ 1 - t (   \lambda_{i,i}^1+\gamma_1^{i,i} )  ]$.

Thus, from the computation in Example~\ref{cplxcoh},
the general $2$-cocycle satisfying the conditions in Theorem~\ref{deformRthm}
Case (ii)
has variables
$( \lambda, \lambda',   \lambda_{i,i}^1,  \lambda_{i,i}^i,
\gamma_1^{1,1}, \gamma_1^{i,i}, \gamma_i^{1,1}, \gamma_i^{i,i}   )$,
with  Eqns. (\ref{Ccocy2eq3}), (\ref{Ccocy2eq4}), (\ref{CR1}),
and (\ref{CR2}).
Equations (\ref{CR1}) and (\ref{CR2}) reduce with Eqn.~(\ref{Ccocy2eq3})
to the same equation
\begin{eqnarray}
3 \gamma_i^{1,1} + 2 \lambda_{i,i}^i + \gamma_i^{i,i}=0, \label{CR4}
\end{eqnarray}
so the deformed $R$ matrix in this case has $5$ free variables.
} \end{example}

\begin{example} \label{gpRex}
{\rm
For $X=k\Z_2$ in Example~\ref{gpex},
the general solutions for the  $2$-cocycles $\phi_1$ and $\phi_2$
with $d^{2,1} = d^{2,2}_{(1)} =d^{2,2}_{(2)}=d^{2,3} =0$ found in
Example~\ref{gpcoh} satisfy the condition (ii)
in Theorem~\ref{deformRthm}: $\phi_1 \tau=\phi_1$, $\tau \phi_2 = \phi_2$.
We check the condition
$(\mu+t \phi_1)(\Delta + t \phi_2)=\delta_1 \id$ on $ \hat{X}$
for some $\delta_1 \in   k[t]/ (t^2)$.
One computes:
\begin{eqnarray*}
(\mu+t \phi_1)(\Delta + t \phi_2)(1)&=& 2+ t\ [  \ (q+
\gamma_1^{x,x} + \lambda' +\lambda_{x,x}^1)  + ( 2r + \lambda +\lambda_{x,x}^x) x \ ] ,
\\  (\mu+t \phi_1)(\Delta + t \phi_2)(x)&=&
2x+ t\ [  \ (r +\gamma_x^{x,x} +2 \lambda )
+ ( 2q + 2 \lambda' ) x \ ] .  \end{eqnarray*}
Thus the general $2$-cocycles  satisfy
$(\mu+t \phi_1)(\Delta + t \phi_2)=\delta_1 \id$,
if and only if the above two values are
multiples of $1$ and $x$, respectively by the same element $\delta_1 \in  k[t]/ (t^2)$.
This condition is written as
\begin{eqnarray*}
2r + \lambda +\lambda_{x,x}^x &=& 0 ,\\
r + \gamma_x^{x,x} +2 \lambda  &=& 0 ,\\ 
\gamma_1^{x,x} + \lambda' +\lambda_{x,x}^1 &=& q + 2 \lambda' .
\end{eqnarray*}
Using Eqns.~(\ref{z2eqn1}) and (\ref{z2eqn2}), one obtains
$ \lambda^x_{x,x}=3\lambda +2\gamma_x^{x,x}$ and one can compute
$\delta_1=2[ 1 + t (\lambda^1_{x,x}+ \gamma_1^{x,x}) ]$. Since the equations listed for Example~\ref{gpcoh} are the only
equations we need for this example as well, the deformed $R$ matrix
has $5$ 
free variables. }\end{example}

\section{Conclusion}

This paper contains a study of the cohomology of Frobenius algebras initiated from the point of view of deformation theory as explicated by diagrammatic techniques. One reason for developing this theory diagrammatically is that we
obtain cocycle deformations of $R$-matrices,
and 
anticipate topological applications.

Several problems remain. Clearly computations of cohomology for
Frobenius algebras formed as matrix algebras would be especially
interesting if these could be done in a fashion that encompassed
all dimensions. The relationship 
between this Frobenius cohomology theory and the
cohomology of the adjoint map in a Hopf algebra deserves study.
Finally, an understanding of the possible knot invariants defined
by YBE solutions defined by skein relations among maps in
Frobenius algebras would be very interesting.

\end{document}